\documentclass[12pt]{article}
\usepackage[utf8x]{inputenc}
\usepackage{amssymb}
\usepackage{color}

\def\emptyset{\varnothing}

\def\1{\mathbf{1}}

\newif\ifhide
\hidetrue

\newtheorem{lem} {Lemma}
\newtheorem{thm} {Theorem}

\newtheorem{rem} {Remark}
\newtheorem{cor} {Corollary}

\newif\ifmoment
\momentfalse

\title{{\normalsize\tt\hfill\jobname.tex}\\
Simple proof of Dynkin's formula
for single--server systems
and polynomial convergence rates }
\author{A. Yu. Veretennikov\footnote{School of
Mathematics, University of Leeds, LS2 9JT, Leeds, UK, \& Institute
of Information Transmission Problems, B. Karetny 19, 127994,
Moscow, Russia,
email: a.veretennikov @ leeds.ac.uk}\hspace{1mm}\footnote{The first author's work was partially supported by RFBR grant 13-01-12447 ofi\_m2.},
\, G. A.
Zverkina\footnote{Department of Applied Mathematics, Moscow State
University of Railway Engineering (MIIT), Obraztsova 9, build. 9,
127994, GSP-4, Moscow, Russia, email: zverkina @ gmail.com} }
\begin{document}

\maketitle

\begin{abstract}
An elementary rigorous justification of Dynkin's identity with an extended generator
based on the idea of a complete probability formula
is given for
queueing systems with a single server and discontinuous intensities of arrivals and serving.
This formula is applied
to the analysis of ergodicity and, in particular, to polynomial
bounds of convergence rate to stationary distribution.

\end{abstract}

\section{Introduction}
In the last decade, systems generalising $M/G/1/\infty,$ or
simply $M/G/1$ (cf. \cite{GK}) -- one of the oldest and most
popular queueing systems -- attracted much attention, see
\cite{Asmussen} -- \cite{fakinos1987},
 \cite{Thor83}
In this paper a single--server system similar to \cite{Ve2013_ait}
is considered,
in which {\em intensities} of new arrivals as well as of their
serving could depend on the ``whole state'' of the system and the
whole state includes the number of customers in the system --
waiting and on service -- {\em and} the elapsed time of the last
serving. Batch arrivals are not allowed. Emphasize that the {\em
elapsed} serving time is assumed known at any moment, but {\em
not} remaining service times for each customer. For definiteness,
the discipline of serving may be regarded as FIFO, although, this
is of no importance for the questions addressed in this paper: all
results are valid for any serving discipline without time
sharing. Recall that despite a non-markovian character of the
model, this model -- as any other process -- may be made Markov by
extending the state space, although not in all cases it may be
helpful; in our case it will be crucial. {\em Dynkin's formula},
or {\em Dynkin's identity} \cite[Ch. 1, \S 3]{Dynkin},
\begin{equation}\label{dynkin1}
E_Xf(X_t) - f(X) = E_X\int\limits_0^t {\mathcal G}f(X_s)\,ds
\end{equation}
for $f$ from some appropriate class of functions and its
non-homogeneous counterpart
\begin{equation}\label{dynkin_t}
E_X\varphi(t,X_t) - \varphi(0,X) = E_X\int\limits_0^t
\left(\frac{\partial}{\partial s}\varphi(s, X_s) + {\mathcal
G}\varphi(s,X_s)\right)\,ds,
\end{equation}
for appropriate functions of two variables $(\varphi(t,X))$ play
a very important role in the analysis of Markov models. Here
${\mathcal G}$ is understood as a {\em generator} or {\em
infinitesimal operator} of the process, and $X$ is a (non-random)
initial value of the process. Note that the formula
(\ref{dynkin_t}) may be considered as a version of (\ref{dynkin1})
written for the couple $(t,X_t)$. In this paper it is shown that
both formulae (\ref{dynkin1})--(\ref{dynkin_t}) hold true for a
large class of functions $f$, $\varphi$ with an {\em
extended generator} ${\mathcal G}$ (see below) and under the
minimal assumptions on regularity of intensities (i.e.,
without any regularity). For short, in the title it was called
Dynkin's identity, but, in fact, Dynkin's formula with
generalised generators will be established. As an application, convergence rate bounds to
stationary distribution are established under appropriate
recurrence conditions similar to those in \cite{Ve2013_ait}, but
without regularity assumed therein.
Minimal or no regularity may
be important, in particular, in queueing systems with control,
where {\em optimal control} is usually discontinuous.

This model may be considered as a partial case of piecewise deterministic Markov processes (PDMP) considered in \cite{Davis}, \cite{GK}, \cite[Ch. 7]{Jacobsen}, \cite{Kalash}, et al. Notice that the method in \cite{Davis},\cite{Jacobsen} is based on martingale theory and requires some further references. In this paper we aim to show that the ``naive'' approach based on a properly designed version of a complete probability formula works well, at least, in our particular model and may be presented in a self--contained way. Also notice that our model allows a presentation as a picewise linear Markov process (PLMP) as in \cite{GK}, \cite{Kalash}, while \cite{Davis} and \cite{Jacobsen} study a more general deterministic behaviour between ``Markov'' jumps.

The standard proof of the formulae
(\ref{dynkin1})--(\ref{dynkin_t}) given in \cite{Dynkin} in ``good cases'' is practically trivial; however,
this triviality is based on the notion of {\em generator} (not extended one), which does
assume some regularity. For example, in the theory of diffusion
processes a similar notion requires continuous coefficients of
drift and diffusion; at the same time, it is well-known that
Dynkin's formula itself still holds true for diffusion processes
with a weaker notion of extended
generator. This is possible due to PDE theory results, namely, {\em Alexandrov'x inequality} and {\em Krylov's estimate}.
In queueing theory there are no similar resources and, as a
consequence, in some classical works a strict justification of
Dynkin's formulae is just dropped.
In the authors' view, this issue must be addressed and the goal of
this paper is to show that such a rigorous justification is well possible practically without -- or with a minimum of -- martingale theory. As a first application, convegence rate bounds to the
stationary regime is shown for a class of models generalising the
system $M/G/1$ similar to \cite{Ve2013_ait} but for a wider class
of intensities. This may be useful in some other areas such as
stochastic averaging and Poisson equations, see, for example,
\cite{kulik2012}. The bounds in \cite{Ve2013_ait} were proposed as
a complement to the results from \cite{Seva}, \cite{Ve77}, et al.
about convergence to stationarity in non-Markov queueing models,
in particulare, of Erlang--Sevastyanov's type. The authors hope to
apply this approach to such systems as well.

The paper consists of the Section 1 -- Introduction, of the setting
and main results in the Section 2, of a set of auxiliary results in the
Section 3, of the proof of Dynkin's formulae in the Section 4,
of the proof of some extension in the Section 5 and of a brief
reminder of the techniques leading to convergence
rate bounds in the Section~6.

\section{The setting and main results}

\subsection{Dynkin's formulae and generalised generator}
Recall \cite{Dynkin} that the generator of a Markov process $(X_t,
\, t\ge 0)$ is an operator \( {\mathcal G}, \) such that for a
sufficiently large class of functions $f$
\begin{equation}\label{dynkin2}
 \sup_X \lim_{t\to 0} \left\|\frac{E_Xf(X_t) - f(X)}{t} -
{\mathcal G}f(X)\right\| = 0
\end{equation}
in the norm of the state space of the process. The notion of
generator does depend on this norm. However, in this paper the
classical generator is not used, so we skip this point here.
Recall further that the class of functions for which
(\ref{dynkin2}) holds true is called {\em domain} of a generator
(of course, in general, it depends on the norm) and in most cases its full
description is a hard task. Fortunately, often it suffices to have
a wide enough class of such functions. On the other hand, an
operator ${\mathcal G}$ is called {\em extended generator} if (\ref{dynkin2}) is replaced by its corollary
(\ref{dynkin1}), also for a wide enough class of functions $f$. It
turns out that in many situations this notion suffices.

Let us present the class of models under investigation in this
paper. The state space is a union of subspaces,
\[
{\mathcal X} = \{(0,0)\} \bigcup \{(n,x): \; n = 1,2,\ldots, N, \;
x\ge 0\}
\]
($1\le N\le \infty$) with topology arising from the metric
($X=(n,x)$, $X'=(n',x')$)
\[
\mbox{dist}(X, X')= |n-n'| + |x-x'|.
\]
Functions of
class $C^1({\mathcal X})$ are understood as functions with
classical continuous derivatives with respect to the variable $x$.
Functions with compact support on ${\mathcal X}$ are understood
as functions vanishing outside some domain bounded in this metric:
for example, $C^1_0({\mathcal X})$ stands for the class of
functions with compact support and one continuous derivative.
There is a generalised Poisson arrival flow with intensity $
\lambda(X), $ where \( X = (n,x) \; \mbox{for $n>0$}, \;
\mbox{and} \; X=(0,0) \; \mbox{for $n=0$}\); group arrival are not allowed.
 If $n>0$, then the
server is serving one customer while all others are waiting in a
queue; if $n=0$ then the server remains idle until the next
arrival and the intensity of such arrival at the state $(0,0)$ is
constant. Here $n$ denotes the total number of customers in the
system and $x$ stands for the elapsed time of the current serving.
Denote $n_t=n(X_t)$ -- the number of customers corresponding to
the state $X_t$ and $x_t=x(X_t)$, the second component of the
process $(X_t)$, i.e. the elapsed time of the current serving. If
$n=0$, then there is no serving. Below, we will investigate formally only the
case with infinite number of waiting places, but the situation where this number
is bounded may be included by assuming $\lambda(n,x)\equiv 0$ for any $n$ large enough.
For
any $X=(n,x)$, intensity of serving $h(X) \equiv h(n,x)$ is
defined; it is also convenient to assume $h(X)=0$ for $n(X)=0$.
Both intensities $\lambda$ and $h$ are understood in the following
way, which is a {\em definition:} on any nonrandom interval of
time $[t,t+\Delta)$, conditional probability given $X_t$ that the
current serving will {\em not} be finished and there will be no
new arrivals reads,
\begin{equation}\label{intu1}
\exp\left(-\int\limits_0^{\Delta} (\lambda+h)
(n_{t},x_{t}+s)\,ds\right).
\end{equation}
This is similar to the approach developed in \cite{Davis}; the main difference is that in this paper the basis for the proof of Dynkin's identity is a complete probability (or complete expectation) formula rather than martingale problem.
In the sequel, $\lambda$ and $h$ are assumed to be {\em bounded}.
In this case, for $\Delta>0$ small enough, the expression in
(\ref{intu1}) may be rewritten as
\begin{equation}\label{intu2}
1-\int\limits_0^{\Delta} (\lambda+h)(n_{t},x_{t}+s)\,ds + O(\Delta^2),
\qquad \Delta\to 0,
\end{equation}
and this what is ``usually'' replaced by
\[
1 - (\lambda(X_t)+h(X_t))\Delta + O(\Delta^2).
\]
In our situation, the latter replacement may be incorrect because
of discontinuities of the functions $\lambda$ and $h$. Emphasize
that from time $t$ until the next jump, the evolution of the
process $X$ is {\em deterministic}. The (conditional given $X_t$)
density of the moment of a new arrival {\em or} of the end of the
current serving after $t$ at $x_t+ z$, $z\ge 0$ equals,
\begin{equation}\label{intu33}
 (\lambda(n_t, x_t+z)+h(n_t, x_t+z))\exp\left(-\int\limits_0^{z}
(\lambda+h)(n_t, x_{t}+s))\,ds\right).
\end{equation}
Further, given $X_t$, the moments of the next ``candidates'' for
jumps up and down are conditionally independent and have the
(conditional -- given $X_t$) density, respectively,
\begin{equation}\label{intu5}
\begin{array}{c}
 \lambda(n_t, x_t+z)\exp\left(-\int\limits_0^{z}
\lambda(n_t, x_t+s)\,ds\right) \; \\
 \\
 \mbox{and} \; \\
 \\
 h(n_t, x_t+z)\exp\left(-\int\limits_0^{z} h(n_t, x_t+s)\,ds\right), \; z\ge 0.\\

\end{array}
\end{equation}
Notice that (\ref{intu33}) does correspond to conditionally
independent densities given in (\ref{intu5}).

For modelling the evolution, both candidate moments should be
realised and the minimal of the two chosen, which will determine
whether the next jump will be, indeed, up or down. The component
$n_t$ is, generally speaking, not Markov, however, the couple
$(X_t = (n_t,x_t))$ in the state space ${\mathcal X}$ {\em is}
Markov by construction and we stress out that continuity of
$\lambda$ or $h$ was not used. The distribution of the process in
the space of trajectories is uniquely determined, that is, if any
other construction is suggested, yet the distributions of all
moments of next jumps up and down ought to be the same.

A natural candidate to the role of the extended generator
is (see a short comment in the beginning of the proof of the
Theorem \ref{thm1})
\begin{eqnarray*}
{\mathcal G}f(X) := \frac{\partial}{\partial x}f(X) + \lambda(X)
(f(X^+) - f(X))
+ h(X) (f(X^-) - f(X)),
\end{eqnarray*}
(recall that $h(0,0) = 0$), where for any $X=(n,x)$,
\[
X^+ := (n+1,x), \quad X^-:= ((n-1)\vee 0,0)
\]
(here $a\vee b = \max(a,b)$).

\begin{thm}\label{thm1}
If the functions $\lambda$ and $h$ are Borel measurable and
bounded, then the formulae (\ref{dynkin1}) and (\ref{dynkin_t})
hold true for any $t>0$ for every $f\in C^1_b({\mathcal X})$ and
$\varphi\in C^1_b([0,\infty)\times {\mathcal X})$, respectively.
Moreover, the process $(X_t, \, t\ge 0)$ is strong Markov with
respect to the filtration \(({\mathcal F}^X_t, \, t\ge 0)\).
\end{thm}

~

Let
$$
L_m(X) = (n+1+x)^m, \quad L_{k,m}(t,X) = (1+t)^k L_{m}(X).
$$
The following extensions of Dynkin's formulae for unbounded
functions hold true: it will be needed in the proof of the Theorem
\ref{thm2}, even though the main body of this proof will be
hidden.
\begin{cor}\label{cor1}
Under the assumptions of the Theorem \ref{thm1},
\begin{eqnarray}\label{M2}
L_{m}(X_t) - L_{m}(X) = \int\limits_0^t \lambda(X_s)
\left[ \left(L_{m}(X^{+}_s) -
L_{m}(X_s)\right) \phantom{\frac{\partial}{\partial} }\right.
 \nonumber \\ \\ \nonumber
\left. + h(X_s) \left(L_{m}(X^{-}_s) - L_{m} (X_s)\right) +
\frac{\partial}{\partial x}L_{m}(X_s)\right]\,ds +M_t,\hspace{0.5cm}
 \end{eqnarray}
with some martingale $M_t$, and also
\begin{eqnarray}\label{M2t}
L_{k,m}(t,X_t) - L_{k,m}(0,X) = \int\limits_0^t \left[\lambda(X_s)
\left(L_{k,m}(s,X^{+}_s) - L_{k,m}(s,X_s)\right) \phantom{\frac{1}{1}}\right.
 \nonumber \\ \\ \nonumber \left.
+ h(X_s) \left(L_{k,m}(s,X^{-}_s) - L_{k,m}(s,X_s)\right) +
\left(\frac{\partial}{\partial x} + \frac{\partial}{\partial
s}\right) L_{k,m}(s,X_s)\right]\,ds +\tilde M_t,\hspace{0.5cm}
 \end{eqnarray}
with some martingale
$\tilde M_t$.
\end{cor}
About martingale language in queueing models see, for example,
\cite{Jacobsen} or \cite{LSh}.

\subsection{Stability and convergence rate}
Denote
\[
\Lambda:= \sup_{n,x: \,n>0}\lambda(n,x).
\]
Recall that the process has no explosion with probability one due
to the boundedness of both intensities, i.e., the trajectory may
have only finitely many jumps on any finite interval of time.
For establishing convergence rate to the stationary regime, we
assume similarly to \cite{Ve2013_ait},
\begin{equation}\label{eq2}
\inf_{n>0} h(n,x) \ge \frac{C_0}{1+x}, \quad x\ge 0.
\end{equation}
\begin{thm}\label{thm2}
Let the functions $\lambda$ and $h$ be Borel measurable and
bounded, the assumption (\ref{eq2}) holds and $C_0$ satisfies the
condition
\begin{equation}\label{condi}
C_0 > 4 (1 + 2 \Lambda),
\end{equation}
and let let the value $k>1$ satisfy
\begin{equation}\label{condi2}
C_0 > 2^{k+1}(1 + \Lambda 2^{k}).
\end{equation}
Then there exists a unique stationary measure $\mu$ and there
exist $m>k$, $C>0$ such that for any $t\ge 0$,
\begin{equation}\label{est}
 \|\mu^{n,x}_t - \mu \|_{TV} \le
 C \,\frac{(1+n+x)^m}{(1+t)^{k+1}},
\end{equation}
where $\mu^{n,x}_t$ is a marginal distribution of the process
$(X_t, \, t\ge 0)$ with the initial data $X=(n,x)\in {\mathcal
X}$.
\end{thm}

\begin{rem}
Under the assumption (\ref{condi}), the inequality (\ref{condi2})
holds true, at least, for all values of $k>1$ close enough to one.
\end{rem}

\begin{rem}
Existence and uniqueness of invariant measure may be established
under weaker assumptions. However, we do not discuss this issue in
this paper, as our goal -- beside Dynkin's formula -- is only convergence rate, for
which we need the condition (\ref{condi}) anyway.
\end{rem}

\section{Auxiliary results}
\noindent The first Lemma is an $L_1$--version of the Theorem
2.8.2 from \cite{Krylov-eng} from integration theory, originally
formulated in $L_2$ as a tool in a Doob's construction related
to {\em stochastic} integrals (the latter are not used here). In
this paper, the Lemma \ref{lem3} will serve as a bridge from the
standard complete probability formula to its integral analogue.
\begin{lem}\label{lem3}
For any $g\in L_1[0,T]$,
\[
\int\limits_0^1 \int\limits_0^T |g_{ \kappa^a_{m'}(s)} - g_s|\,ds\,da \to 0,
\]
and for almost all $a\in [0,1]$,
\[
 \int\limits_0^T |g_{ \kappa^a_{m'}(s)} - g_s|\,ds \to 0,
\]
over some subsequence $m'\to\infty$, where it was denoted
\[
 \kappa^a_m(s):= [2^m(s+a)]2^{-m}-a,
\]
where, in turn, $[a]$ is the integer value of $a\in R$ (i.e., the
nearest integer approximation of $a$ from the left). It is
accepted that outside $[0,T]$ the function $g$ equals identically
zero.

\end{lem}
{\em The proof} is an exact repetition of the calculus in the
proof of \cite[Theorem 2.8.2]{Krylov-eng} with a tiny difference
of $L_2$ replaced by $L_1$, which does admit this small change.
Hence, the details are dropped. For convenience of the reader we
only recall that it follows from approximation of $g\in L_1$ by
continuous functions, which are dense in $L_1$ and establishing
the statement for continuous functions, which is evident. In the
proof of the Theorem \ref{thm1} below a similar trick will be used
for establishing a similar auxiliary assertion (\ref{la4}), in
which a one-dimensional nature of time seems to be important.

~

\noindent The second Lemma provides a rigorous proof of the
well-known properties that probability of ``one event'' on a small
nonrandom interval of length $\Delta$ is of the order $O(\Delta)$
and probability of ``two or more events'' on the same interval is
of the order $O(\Delta^2)$; although this is a ``common
knowledge'' in queueing theory, the authors believe that for
discontinuous intensities it must be justified.

In the next Lemma we accept the following convention: all trajectories of the process are right continuous with left limits; respectively, if the value of the process $X_t$ at $t$ is given, the phrase ``no jumps on $(t,t+s]$ relates to no jumps after the moment $t$, but not at $t$, where we are not aware whether or not $X_{t-}=X_t$.  At the same time, ``no jumps on $(0,s]$'' includes no jump at zero, as negative values of $t$ are not allowed, so that automatically $X_0=X_{0+}$.

\begin{lem}\label{lem1}
Under the assumptions of the Theorem \ref{thm1}, for any $t\ge 0$,
\begin{eqnarray}
P_{X_{t}}(\mbox{no jumps on $(t, t+\Delta]$}) =
\exp\left(-\int\limits_0^\Delta (\lambda+h)(X_t+s)\,ds\right) = 1 + O(\Delta),\hspace{0.6cm}\label{z0}\\
P_{X_{t}}(\mbox{at least one jump on $(t, t+\Delta]$}) = O(\Delta),\hspace{6cm}\label{z1}\\
P_{X_{t}}(\mbox{exactly one jump up \& no down on $(t,
t+\Delta]$}) = \int\limits_0^\Delta \lambda(X_t+s)\,ds + O(\Delta^2),\label{z1up}\\
P_{X_{t}}(\mbox{exactly one jump down \& no up on $(t,
t+\Delta]$}) = \int\limits_0^\Delta h(X_t+s)\,ds + O(\Delta^2),\hspace{0.1cm}\label{z1down}
\end{eqnarray}
and
\begin{eqnarray}
P_{X_{t}}(\mbox{at least two jumps on $(t, t+\Delta]$}) =
O(\Delta^2).\hspace{5.7cm}\label{z2}
\end{eqnarray}
In all cases above, $O(\Delta)$ and $O(\Delta^2)$ are uniform
with respect to $X_{t}$ and only depend on the norm
$\sup_{X}(\lambda(X)+h(X))$, that is, there exist $C>0, \,
\Delta_0>0$ such that for any $X$ and any $\Delta<\Delta_0$,
\begin{eqnarray}\label{z_uni}
 \limsup_{\Delta\to 0} \left\{\Delta^{-1}P_{X}(\mbox{at least one
jump on $(0,\Delta]$})\phantom{\int\limits_1^1}\right.\hspace{2cm}
 \nonumber \\ \nonumber \\ \nonumber
+ \Delta^{-2} P_{X}(\mbox{at least two jumps on $(0, \Delta]$})\hspace{2.5cm}
  \nonumber \\  \nonumber \\\nonumber
+ \Delta^{-2}\left[P_{X_{t}}(\mbox{one jump up \& no down on $(t,
t+\Delta]$}) - \int\limits_0^\Delta \lambda(X_t+s)\,ds \right]
 \nonumber  \\ \nonumber \\
\left.+ \Delta^{-2}\left[P_{X_{t}}(\mbox{one jump down \& no up on
$(t, t+\Delta]$}) - \int\limits_0^\Delta h(X_t+s)\,ds \right]\right\} \nonumber \\ \nonumber\\
\le C< \infty.\hspace{5cm}
\end{eqnarray}
\end{lem}
{\bf Proof}. {\bf A.} The assertions (\ref{z0}), as well as
(\ref{z1}), follows from the construction (see (\ref{intu1}) and
(\ref{intu2}), since $\lambda$ and $h$ are bounded. (It also
follows from the next steps of the proof.)

~

\noindent {\bf B.} Now it is convenient to proceed with
(\ref{z2}); it will be used in the proof of (\ref{z1up}) and
(\ref{z1down}). It suffices to consider $t=0$ and any fixed
initial value $X$. Denote
\[
s^N_j = \frac{j\,\Delta}{2^N}, \qquad 1\le j\le 2^N, \;\; N\ge 1,
\]
and let
\begin{eqnarray*}
B :=\{\mbox{two or more jumps -- up or down -- on $(0,\Delta]$}\},\hspace{2cm}\\
\\
B^N:=\bigcup_{j} \left\{\mbox{no jumps on $(0,s^N_{j-1}]$,}\right.
 \mbox {at
least one jump on $(s^N_{j-1},s^N_j]$}
 \\ \\
\left. \mbox{and
 at least one jump on $(s^N_j, \Delta]$}\right\},\hspace{3.2cm}
\end{eqnarray*}
and
\[
 p:= P_{X}(B), \qquad p_N := P_{X}(B^N).\hspace{5cm}
\]
Then,
\begin{eqnarray*}
p_N= \sum_{j\le 2^N} P_X(\mbox{no jumps on $(0,s^N_{j-1}]$,} \mbox{at
least one jump on $(s^N_{j-1},s^N_j]$}
 \\ \\
\mbox{and at least one jump on $(s^N_j, \Delta]$}).\hspace{3.2cm}
\end{eqnarray*}
Notice that for any nonrandom $t$, $P(X_t\not= X_{t-}) =0$. We
have,
\[
B^N \subset B,
\]
and, moreover,
\[
B \setminus B^N \downarrow \emptyset, \quad N\uparrow \infty.
\]
The latter holds because for any two jumps there exists $N$ such
that they would be covered by two different intervals, say,
$(s^N_{j-1},s^N_j]$ and $(s^N_{k-1},s^N_k]$ with $k\not=j$. In
other words, for each $\omega\in B$ there exists $N$ such that
$\omega \in B^N$. Hence, by continuity of a probability measure,
\[
 p_N\to p, \quad N\to\infty.
\]
From the Markov property and (\ref{z0}--\ref{z1}), if $\Delta$ is
small enough then we have,
\begin{eqnarray*}
 p_N
= \sum_{j \le 2^N} E_X E \left( 1(\mbox{no jumps
on $(0,s^N_{j-1}]$}) \right. \hspace{6.6cm}
\\\\
\times1(\mbox{at least one jump on
$(s^N_{j-1},s^N_j]$} \mid {\mathcal F}^{X}_{s^N_j})
\times \left.1 (\mbox{at least one jump on $(s^N_j,
\Delta]$})\right)
 \\\\
= \sum_{j \le 2^N} E_X \left(1(\mbox{no jumps on $(0,s^N_{j-1}]$})\right.  \hspace{7.5cm}
\\\\
\times \left.1(\mbox{at least one jump on $(s^N_{j-1},s^N_j]$})\right)
\times E\left(1 (\mbox{at least one jump on $(s^N_j,
\Delta]$}\mid X_{s^N_j}\right)
 \\\\
= \sum_{j \le 2^N} E_X \left(1(\mbox{no jumps on $(0,s^N_{j-1}]$}) \right.  \hspace{7.5cm}
\\\\
\times\left. 1(\mbox{at least one jump on $(s^N_{j-1},s^N_j]$}) \right) \times
\left(1-\exp\left(-\int\limits_{s^N_j}^{\Delta}(\lambda+h)
(X_{s^N_j}+s)\,ds\right)\right)
 \\\\
\le C \Delta \, \sum_{j \le 2^N} E_X 1\left( \mbox{no jumps on
$(0,s^N_{j-1}]$}, \mbox{at least one jump on
$(s^N_{j-1},s^N_j]$}\right)
 \\\\
= C \Delta \,P_X (\mbox{at least one jump up on
$(0,s^N_N]$})\hspace{6.5cm}
 \\\\
= C \Delta \,\left(1-\exp\left(-\int\limits_0^{\Delta}(\lambda+h)
(X_{}+s)\,ds\right)\right)
\le C^2 \Delta^2.  \hspace{4.3cm}
\end{eqnarray*}
Important is that here the last value of $C^2$ is the same for all
values of $N$ and $X$. Hence, by monotone convergence, we also
have
\[
 p \le C^2\Delta^2,
\]
as required. Notice that (\ref{z1up}) and (\ref{z1down}) were not
used so far.

~

\noindent {\bf C.} Let us show ({\ref{z1up}). Informally,
\begin{eqnarray*}
P_{X_{t}}(\mbox{exactly one jump up \& no down on $(t,
t+\Delta]$}) \hspace{5cm}
 \\\\
= \int\limits_0^\Delta \lambda(X_t+s)\, \exp\left(-\int\limits_0^s
\lambda(X_t+r)\,dr\right) \exp\left(-\int\limits_s^\Delta \lambda(X_s^+
+r)\,dr\right)\hspace{0.5cm}
 \\\\
\times \exp\left(-\int\limits_0^s h(X_t+r)\,dr\right)
\exp\left(-\int\limits_s^\Delta h(X_s^+ +r)\,dr\right)\,ds
\hspace{0.5cm}
 \\\\
= \int\limits_0^\Delta \lambda(X_t+s)\,ds + O(\Delta^2),
\hspace{7.9cm}
\end{eqnarray*}
as required and with a uniform $O(\Delta^2)$. Some drawback of this
explanation is that it uses an integral identity as if it were
complete probability formula, while the latter formula is stated for a split of $\Omega$ into finitely or countably many disjoint events only. Usually there is no problem with such integration, at least for Riemann integrable functions, since integration is a limit of Darboux sums. However, in our case we assume the integrands, in general, only Lebesque integrable. So, let us do this estimate more
rigorously. We will use the already established formula (\ref{z2}) and
instead of (\ref{z1up}) we will estimate a slightly different probability
of {\em at least one jump up earlier than down} on $(0,\Delta]$
(we mean that jump down may or may not occur on this interval).
According to our construction, in particular, (\ref{intu33}), the
value of this probability may be written precisely by using the
(conditionally independent given $X_t$) densities (\ref{intu5}):
\[
f^u(z):=\lambda(X_t+z)\exp\left(-\int\limits_0^{z}
\lambda(X_{t}+s)\,ds\right)
\]
and
\[
f^d(z):= h(X_t+z)\exp\left(-\int\limits_0^{z}
h(X_{t}+s)\,ds\right).
\]
We have,
\begin{eqnarray*}
 P_X(\mbox{at least one jump up earlier than down on $(0,\Delta]$})\hspace{3cm}
 \\\\
= \int\limits_0^\infty \int\limits_0^\infty 1(z^1<z^2)1(z^1\le \Delta)
f^u(z^1)f^d(z^2)\,dz^1dz^2\hspace{3cm}
 \\\\
\le \int\limits_0^\infty \int\limits_0^\infty 1(0<z^2) 1(z^1\le \Delta)
f^u(z^1)f^d(z^2)\,dz^1dz^2\hspace{3.2cm}
 \\\\
= \int\limits_0^\Delta f^u(z^1) \,dz^1 = \int\limits_0^\Delta
\lambda(X+z)\exp\left(-\int\limits_0^{z} \lambda(X+s)\,ds\right) \,dz\hspace{1.0cm}
 \\\\
\le \int\limits_0^\Delta \lambda(X+z) \,dz.\hspace{8.2cm}
\end{eqnarray*}
On the other hand,
\begin{eqnarray*}
 P_X(\mbox{at least one jump up earlier than down on $(0,\Delta]$})\hspace{4cm}
 \\\\
= \int\limits_0^\infty \int\limits_0^\infty 1(z^1<z^2)1(z^1\le \Delta)
f^u(z^1)f^d(z^2)\,dz^1dz^2\hspace{4.1cm}
 \\\\
\ge \int\limits_0^\infty \int\limits_0^\infty 1(\Delta<z^2) 1(z^1\le \Delta)
f^u(z^1)f^d(z^2)\,dz^1dz^2\hspace{4.2cm}
 \\\\
= \int\limits_0^\Delta f^u(z^1) \,dz^1 \times \int\limits_\Delta^\infty f^d(z^2)
\,dz^2\hspace{7.2cm}
 \\\\
\ge \int\limits_0^\Delta \lambda(X+z)\exp\left(-\int\limits_0^{\Delta}
\lambda(X+s)\,ds\right) \,dz \times\exp\left(-\int\limits_0^{\Delta}
h(X+s)\,ds\right)
 \\\\
= \int\limits_0^\Delta \lambda(X+z) \,dz\times
\exp\left(-\int\limits_0^{\Delta} (\lambda + h)(X+s)\,ds\right)\hspace{3.5cm}
 \\\\
= (1+O(\Delta)) \int\limits_0^\Delta \lambda(X+z) \,dz = \int\limits_0^\Delta
\lambda(X+z) \,dz + O(\Delta^2),\hspace{2.4cm}
\end{eqnarray*}
with uniform $O(\Delta)$ and $O(\Delta^2)$, as required. Since we
already know that probability of two or more jumps is of the order
$\Delta^2$ -- see (\ref{z2}) -- this justifies (\ref{z1up}).

~

\noindent {\bf D.} The statement ({\ref{z1down}) is established similarly.
Note that in the previous steps, expressions like
$\displaystyle 1-\exp\left(-\int\limits_0^\Delta \psi(t)\,dt\right)$ with some non-negative bounded $\psi(t)$ are evaluated, say, $\psi(t)<M$, $\forall t$. In this case, this exponent may be expanded via  Taylor's series of the variable $\displaystyle \int\limits_0^\Delta \psi(t)\,dt$. For $\Delta <2/M$ this series satisfies the well-known condition for the alternating series test, which implies strict bounds $\displaystyle 1-\exp\left(-\int\limits_0^\Delta \psi(t)\,dt\right)<M\Delta$ and $\displaystyle 0<\left(1-\exp\left(-\int\limits_0^\Delta \psi(t)\,dt\right)\right)-\int\limits_0^\Delta \psi(t)\,dt< M^2\Delta^2/ 2$. This shows that the estimates in the proof above and in
({\ref{z_uni}}) are uniform, as required.
The Lemma \ref{lem1} is proved.

~

\begin{lem}\label{Le2}
Under the assumptions of the Theorem \ref{thm1}, for any $f\in C_b({\mathcal X})$ the function
\linebreak $T_tf(X) = E_Xf(X_t)$ is continuous in $t$.
\end{lem}
{\bf Proof}. It suffices for any $f\in C_b({\mathcal X})$ and for
small $\Delta>0$ to show that for any $t\ge 0$,
\begin{equation}\label{Feller}
 |E_X (f(X_{t+\Delta}) - f(X_{t}))| \to 0, \quad
\Delta \to 0.
\end{equation}
and for any $t>0$,
\begin{equation}\label{Feller2}
|E_X (f(X_{t-\Delta}) - f(X_{t})|\to 0, \quad \Delta \to 0.
\end{equation}
The convergence (\ref{Feller}) follows straightforward from the
first statement of the Lemma \ref{lem1} by virtue of Lebesgue's
dominated convergence theorem, due to the cadlag property of the
process. Similarly, we have,
\[
|E_X (f(X_{t-\Delta}) - f(X_{t-})|\to 0, \quad \Delta \to 0.
\]
However, by construction (see (\ref{intu1})), jump at any nonrandom
time $t$ has probability zero, that is,
\[
P(X_t=X_{t-}) = 1.
\]
This implies (\ref{Feller2}). The Lemma \ref{Le2} is proved.

~

\noindent Denote $X'=(n',x')\uparrow X=(n,x)$ iff $x'\uparrow x$
and $n' = n$ when $x'$ is close enough to $x$. Similarly,
$X'=(n',x')\downarrow X=(n,x)$ iff $x'\downarrow x$ and $n' = n$
for $x'$ close enough to $x$.
\begin{lem}\label{lem4}
Under the assumptions of the Theorem \ref{thm1} the process
$(X_t,\,t\ge 0)$ is Feller, that is, $T_tf(\cdot)\in C_b({\mathcal
X})$ for any $f\in C_b({\mathcal X})$.
\end{lem}
{\em Proof.} \noindent {\bf A.} Denote $\phi(t,X):=E_Xf(X_t)$. Let
$(n', x')=X'\uparrow X=(n,x)$. Without loss of generality we may
assume $n'\equiv n$ and $x'\uparrow x$. Let $s:=x-x'>0$. Note that
for the process with initial data $X'=(n',x')$,
\[
x_s \,
1(\mbox{no jumps on $(0,s]$})
= x 
\,1(\mbox{no jumps
on $(0,s]$}).
\]
So, for $x'$ close enough to $x$, by the Lemma \ref{lem1},
\begin{eqnarray*}
\phi(t,X') = E_{X'}f(X_t) \hspace{8cm}\\\\
=E_{X'}f(X_t)1(\mbox{no jumps on
$(0,s]$})\hspace{3.5cm}
 \\\\
+ E_{X'}f(X_t)1(\mbox{at least one jump on $(0,s]$})\hspace{1cm}
 \\\\
= E_{X'}1(\mbox{no jumps on $(0,s]$})E_{X_s}f(X_{t-s}) + O(s)\hspace{1.1cm}
 \\\\
=E_{X'}1(\mbox{no jumps on $(0,s]$})E_{X}f(X_{t-s}) + O(s)\hspace{1.3cm}
 \\\\
= E_{X'}1(\mbox{no jumps on $(0,s]$})\phi(t-s,X) + O(s)\hspace{1.3cm}
 \\\\
= (1-O(s))\phi(t-s,X) + O(s).\hspace{3.7cm}
\end{eqnarray*}
By virtue of the Lemma \ref{Le2} and due to our notation $s=x-x'$,
\[
\lim_{X'\uparrow X} \phi(t,X') = \lim_{s\downarrow 0} \phi(t-s,X)
= \phi(t,X).
\]

~

\noindent {\bf B.} Now let $X'\downarrow X$. In this case denote
$s:=x'-x$ ($s>0$). Then similarly to the above,
\[
E_{X}f(X_t) = O(s) + (1 - O(s))E_{X'}f(X_{t-s})
\]
Hence, we have,
\[
\phi(t,X) = \phi(t-s, X') + O(s).
\]
Since by the Lemma \ref{Le2},
$$
\phi(t,X') = \phi(t-s, X') + o(1),
$$
we obtain,
\[
\phi(t,X') = \phi(t-s, X') + o(1) = \phi(t, X) + O(s) + o(1),
\]
or, equivalently,
\[
\lim_{X'\downarrow X}\phi(t,X') = \phi(t, X).
\]
This completes the proof of the Lemma \ref{lem4}.

\begin{rem}
Note that in \cite{Davis} Feller's property is proved for a more general model, however, under the additional condition that intensities are continuous. In our case continuity is not necessary -- which is in line with the very idea of construction of a process with discontinuous intensities -- and in the sequel this Feller's property will be used for establishing strong Markov property by using classical tools. On the other hand, in \cite{Davis} strong Markov property is established independently from Feller's one.
\end{rem}

\section{Proof of Theorem \ref{thm1}}
As was already mentioned above, both formulae (\ref{dynkin1}) and
(\ref{dynkin_t}) are, actually, versions of complete probability
formula (it would be better to say, complete expectation). If it
were possible to replace the integrals by series, they both would
have been, indeed, complete probability (expectation) formulae as
suggested. However, in our case, integrals and series are not the
same things. Notice that for {\em continuous} and bounded
$\lambda$ and $h$, the formula (\ref{dynkin1}) with any $f\in
C^1_b({\mathcal X})$ is, indeed, a simple corollary of a standard
analysis of probabilities of transitions over a small time, as in
the derivation of the Kolmogorov forward equations. This is why
${\mathcal G}$ may be regarded as a natural candidate to a
generalised generator for discontinuous intensities, too.

~

\noindent {\bf Ia.} In the first part of the proof we assume $f\in
C^1_0({\mathcal X})$, rather than $f\in C^1_b({\mathcal X})$ in
the formula (\ref{dynkin1}). In some lines it is useful to know
that $f$ and its derivative are uniformly continuous. In the end
of the proof, the assertion (\ref{dynkin1}) will be extended from
$C^1_0({\mathcal X})$ to the whole $C^1_b({\mathcal X})$. Let $\Delta:=2^{-m}t$,
\[
t_0=0, \;\; t_i = \kappa^a_m\left(\frac{it}{2^m}\right) \vee 0\equiv ([it+a 2^m]2^{-m}-a)\vee 0, \; 0\le i\le 2^m, \quad
t_{2^m+1}=t.
\]
Let $\Delta_i:= t_{i+1}-t_i$. Notice that $\Delta_i=\Delta$, $1\le
i\le 2^m$ and $\Delta_i\le \Delta$ for any $i$. \\
Consider the difference
\begin{eqnarray*}
 E_X f(X_t) - f(X) = E_X \sum_{i=0}^{2^m} E_{X_{t_i}}
(f(X_{t_{i+1}})-f(X_{t_{i}})).
\end{eqnarray*}
Consider one term from this sum. Emphasize that they are all
treated similarly. This term may be split into four parts:
\begin{eqnarray*}
 E_{X_{t_i}}(f(X_{t_{i+1}})-f(X_{t_{i}})) = E_{X_{t_i}}
(f(X_{t_{i+1}})-f(X_{t_{i}}))1(\mbox{no jumps on
$(t_{i},t_{i+1}]$})
 \\\\
 + E_{X_{t_i}}(f(X_{t_{i+1}})-f(X_{t_{i}}))1(\mbox{one
jump up \& no down on $(t_{i},t_{i+1}]$})
 \\\\
 + E_{X_{t_i}}(f(X_{t_{i+1}})-f(X))1(\mbox{one jump
down \& no up on $(t_{i},t_{i+1}]$})\hspace{0.2cm}
 \\\\
 + E_{X_{t_i}}(f(X_{t_{i+1}})-f(X_{t_{i}}))1(\mbox{at
least, two jumps on $(t_{i},t_{i+1}]$})\hspace{1cm}
 \\\\
 \equiv I^i_1 + \ldots + I^i_4. \hspace{5cm}
\end{eqnarray*}

\noindent {\bf Ib.} We have, with notation $X+s:= (n,x+s)$ for any $X=(n,x)$,
\begin{eqnarray}\label{i1}
 I^i_1 = \exp\left(-\int\limits_0^{\Delta_i}
(\lambda(X_{t_i}+s)+h(X_{t_i}+s))\,ds\right)
(f(X_{t_i}+\Delta_i)-f(X_{t_i}))
 \nonumber \\ \\ \nonumber
 = (1+O(\Delta))\int\limits_0^{\Delta_i}\frac{\partial}{\partial
x}f(X_{t_i}+s)\,ds. \hspace{3cm}
\end{eqnarray}
Here it is likely that we could use the Lemma \ref{lem3}, but we
prefer easier methods where possible. Notice that
\begin{eqnarray}\label{i11}
 \int\limits_0^{\Delta_i} \frac{\partial}{\partial x}f(X_{t_i}+s)\,ds
 = \int\limits_0^{\Delta_i} \frac{\partial}{\partial
x}f(X_{t_i+s})\,ds
 \nonumber \\ \\ \nonumber
+ \int\limits_0^{\Delta_i} \left(\frac{\partial}{\partial x}f(X_{t_i}+s)
- \frac{\partial}{\partial x}f(X_{t_i+s})\right)\,ds.
\end{eqnarray}
After summation and taking expectation, this term gives us
\begin{eqnarray*}
E_X\sum_{i=0}^{2^m} I^i_1 = E_X \sum_{i=0}^{2^m}
\int\limits_0^{\Delta_i} \frac{\partial}{\partial x}f(X_{t_i}+s)\,ds
 = E_X \sum_{i=0}^{2^m} \int\limits_0^{\Delta_i} \frac{\partial}{\partial
x}f(X_{t_i+s})\,ds
 \nonumber \\ \\ \nonumber
+ E_X \sum_{i=0}^{2^m} \int\limits_0^{\Delta_i}
\left(\frac{\partial}{\partial x}f(X_{t_i}+s) -
\frac{\partial}{\partial x}f(X_{t_i+s})\right)\,ds.
 \hspace{1cm}
\end{eqnarray*}
By virtue of the cadlag property of $(X_s, \, s\ge 0)$ and $f\in
C^1_b$, we have by Lebesgue's dominated convergence Theorem,
\begin{eqnarray*}
E_X \sum_{i=0}^{2^m} \int\limits_0^{\Delta_i}
\left(\frac{\partial}{\partial x}f(X_{t_i}+s) -
\frac{\partial}{\partial x}f(X_{t_i+s})\right)\,ds\hspace{6cm}
 \\\\
= E_X \int\limits_0^t \left(\frac{\partial}{\partial
x}f(X_{\kappa^a_m(s)}+s-\kappa^a_m(s) - \frac{\partial}{\partial
x}f(X_{s})\right)\,ds\hspace{2.8cm}
 \\\\
= E_X \int\limits_0^t \left(\frac{\partial}{\partial
x}f(X_{\kappa^a_m(s)}+s-\kappa^a_m(s) - \frac{\partial}{\partial
x}f(X_{s-})\right)\,ds \to 0, \;\; \Delta \to 0,
\end{eqnarray*}
the latter equality a.s., because
\[
P_X(X_{\kappa^a_m(s)}+s-\kappa^a_m(s)-X_{s-} \to 0, \quad
m\to\infty)=1.
\]
Hence -- and due to $\int\limits_0^t |f(X_s)-f(X_{s-})|\,ds = 0$ a.s. -- we get,
\begin{eqnarray*}
E_X \sum_{i=0}^{2^m} I^i_1 = E_X \int\limits_0^t
\frac{\partial}{\partial x}f(X_{s})\,ds + o(1), \;\; \Delta \to 0.
\end{eqnarray*}

~

\noindent {\bf Ic.} Further,
\begin{eqnarray}\label{Ic}
 I^i_2
= E_{X_{t_i}}(f(X_{t_{i+1}})-f(X_{t_{i}}))1(\mbox{one jump up \&
no down on $(t_{i},t_{i+1}]$})\hspace{1cm}
 \nonumber \\ \nonumber \\ \nonumber
\equiv E_{X_{t_i}}(f(X^+_{t_{i}})-f(X_{t_{i}}))1(\mbox{one jump up \&
no down on $(t_{i},t_{i+1}]$})\hspace{0.5cm}
 \\ \nonumber \\ \nonumber
+ E_{X_{t_i}}(f(X_{t_{i+1}})-f(X^+_{t_{i}}))1(\mbox{one jump up \&
no down on $(t_{i},t_{i+1}]$})
 \\ \nonumber \\ \nonumber
= (f(X^+_{t_i}) - f(X_{t_i})) \left(\int\limits_0^{\Delta_i}
\lambda(X_{t_i}+s) \,ds + O(\Delta^2)\right)\hspace{3cm}
 \nonumber \\\\ \nonumber
+ E_{X_{t_i}}(f(X_{t_{i+1}})-f(X^+_{t_{i}}))1(\mbox{one jump up \&
no down on $(t_{i},t_{i+1}]$}),
\end{eqnarray}
by the Lemma \ref{lem1}. After summation and taking expectation,
the first term here gives
\begin{eqnarray*}
E_X \int\limits_0^t (f(X^+_{\kappa^a_m(s)}) - f(X_{\kappa^a_m(s)}))
\lambda(X_{\kappa^a_m(s)}+s - \kappa^a_m(s)) \,ds + O(\Delta)),
\end{eqnarray*}
and we will explain below -- via the Lemma 1 -- why it approaches
the desired
\begin{eqnarray*}
E_X \int\limits_0^t (f(X^+_{s-}) - f(X_{s-})) \lambda(X_{s}) \,ds = E_X
\int\limits_0^t (f(X^+_{s}) - f(X_{s})) \lambda(X_{s}) \,ds, \;\;
\mbox{a.s.}
\end{eqnarray*}
The absolute value of the second term in (\ref{Ic}) does not
exceed
\[
\rho_f(\Delta) E_{X_{t_i}}1(\mbox{one jump up \& no down on
$(t_{i},t_{i+1}]$}),
\]
where $\rho_f$ is the modulus of continuity of the function $f$.
By virtue of the Lemma \ref{lem1} (see (\ref{z1up})),
after summation and expectation this gives us
\[
\rho_f(\Delta) E_{X}\int\limits_0^t \lambda(X_{\kappa^a_m(s)}+s -
\kappa^a_m(s))\,ds + O(\Delta) = o(1), \quad \Delta\to 0.
\]

~

\noindent {\bf Id.} Similarly
with some $\theta\in (0,\Delta_i)$,
\begin{eqnarray*}
I^i_3 = (1 - \exp\left(-\int\limits_0^{\Delta_i} h(X_{t_i}+s) \,ds\right)
+ O(\Delta^2))
 (f(X^-_{t_i}+\theta) - f(X_{t_i}))
 \\\\
 = \left(\int\limits_0^{\Delta_i} h(X_{t_i}+s) \,ds + O(\Delta^2)\right)
 (f(X^-_{t_i}) - f(X_{t_i}) + \rho_f(\Delta)).
\end{eqnarray*}
Since $\rho_f(\Delta) = o(1)$, we get,
\begin{eqnarray*}
 E_X \sum_{i=0}^{2^m} I^i_3=\hspace{12cm}
 \\\\
= E_X \int\limits_0^t h(X_{\kappa^a_m(s)}+s - \kappa^a_m(s)) \,ds
 (f(X^-_{\kappa^a_m(s)}) - f(X_{\kappa^a_m(s)}) + o(1)).
\end{eqnarray*}
It will be shown below that the main term here approaches
\begin{eqnarray*}
E_X \int\limits_0^t h(X_{s}) (f(X^-_{s-}) - f(X_{s-}))\,ds
 \stackrel{\mbox{\small a.s.}}{=} E_X \int\limits_0^t h(X_{s}) (f(X^-_{s}) - f(X_{s})\,ds.
\end{eqnarray*}

~

\noindent {\bf Ie.} Finally (see the Lemma \ref{lem1}), with
uniform $O(\Delta^2)$,
\[
 I^i_4 = O(\Delta^2) \Longrightarrow E_X \sum_{i=0}^{2^m} I^i_4 = O(\Delta).
\]

~

\noindent {\bf If.} Now in the expression arising from the terms
$E_X I^i_2$, as well as from $E_X I^i_3$, we would like to replace
$f(X^+_{\kappa^a_m(s)})
 - f(X_{\kappa^a_m(s)})$ by
$f(X^+_{s})
 - f(X_{s})$ and
$f(X^-_{\kappa^a_m(s)})
 - f(X_{\kappa^a_m(s)})$ by
$f(X^-_{s})
 - f(X_{s})$,
respectively; also,
$h(X_{\kappa^a_{m'}(s)})$ and $\lambda(X_{\kappa^a_{m'}(s)})$ will be
replaced by
$h(X_{s})$ and $\lambda(X_{s})$, respectively. By virtue of continuity of $f$ and
cadlag property of the process, we estimate
for almost every $a\in [0,1]$,
\begin{eqnarray*}
 \left(E_X \int\limits_0^t \lambda(X_{\kappa^a_m(s)})
(f(X^+_{\kappa^a_m(s)})
 - f(X_{\kappa^a_m(s)}))\,ds\right.
 \\\\
\left. - E_X \int\limits_0^t \lambda(X_{s})(f(X^+_{s})
 - f(X_{s}))\,ds\right)
 \to 0, \quad m \to\infty,
\end{eqnarray*}
and
\begin{eqnarray*}
 \left(E_X \int\limits_0^t h(X_{\kappa^a_{m'}(s)})(f(X^-_{\kappa^a_{m'}(s)})
 - f(X_{\kappa^a_{m'}(s)}))\,ds\right.
 - \\\\
 \left.- E_X \int\limits_0^t h(X_{s})(f(X^-_{s+})
 - f(X_{s}))\,ds\right)
 \to 0, \quad m' \to\infty,
\end{eqnarray*}
where $m'\to \infty$ is some subsequence guaranteed by the Lemma \ref{lem3}.
Indeed, say, for $h$, by the Lemma \ref{lem3} for a.e. $a$,
\begin{eqnarray*}
 E_X \int\limits_0^t |(h(X_{\kappa^a_{m'}(s)}) - h(X_s))
(f(X^+_{\kappa^a_{m'}(s)})
 - f(X_{\kappa^a_{m'}(s)}))|\,ds
 \\\\
 \le C_f \, E_X \int\limits_0^t |(h(X_{\kappa^a_{m'}(s)}) - h(X_s))|\,ds
 \to 0, \quad m' \to\infty.
\end{eqnarray*}
Further, by continuity of $f$ and cadlag property of $X$,
\begin{eqnarray*}
\left( E_X \int\limits_0^t h(X_{s})(f(X^+_{\kappa^a_{m'}(s)})
 - f(X_{\kappa^a_{m'}(s)}))\,ds -\right.\hspace{1cm}
\\\\
\left. - E_X \int\limits_0^t h(X_{s})(f(X^+_{s-})
 - f(X_{s-}))\,ds\right)
 \to 0, \quad m' \to\infty,
\end{eqnarray*}
by virtue of the Lebesgue dominated (bounded) convergence theorem, as
$f(X_{\cdot})$ is cadlag, both $h$ and $f$ are bounded and
\[
 \kappa^a_{m'}(s) \to s-, \quad m'\to\infty.
\]
Finelly, since \(P(X_s\not=X_{s-})=0, \, \forall s\),
we get,
\begin{eqnarray*}
 E_X \left|\int\limits_0^t h(X_{s})(f(X^+_{s-})
 - f(X_{s-}))\,ds - \int\limits_0^t h(X_{s})(f(X^+_{s})
 - f(X_{s}))\,ds
\right| = 0.
\end{eqnarray*}

~

\noindent {\bf Ig.} Further, we would like to have $ \lambda(X_s)$
under the integral instead of $\lambda(X_{t_i}+(s-t_i))$ in each term (\ref{Ic})
\[
(f(X^+_{t_i}) - f(X_{t_i})) \left(\int\limits_0^{\Delta_i}
\lambda(X_{t_i}+s) \,ds + O(\Delta^2)\right).
\]
So, let us estimate the difference
\begin{eqnarray*}
 E_X\sum_{i=0}^{2^m} \left(\int\limits_{t_i}^{t_{i+1}}
\lambda(X_{t_i}+(s-t_i)) \,ds - \int\limits_{t_i}^{t_{i+1}}
 \lambda(X_s)\,ds\right)
 (f(X^+_{t_i}) - f(X_{t_i})).
\end{eqnarray*}
Recall that $t_i = \kappa^a_m\left(\frac{it}{2^m}\right)$. Let us show that
for almost any $a\in [0,1]$,
\begin{eqnarray}\label{main_dif1}
 E_X \sum_{i=0}^{2^m} \left|\int\limits_{t_{i}}^{t_{i+1}} \lambda(X_{t_i}+(s-
t_i)) \,ds
 - \int\limits_{t_{i}}^{t_{i+1}} \lambda(X_s)\,ds\right|\,
 |f(X^+_{t_{i}}) - f(X_{t_{i}})|= \hspace{2.1cm}
 \nonumber \\\\\nonumber
E_X\sum_{i=0}^{2^m} \left|\int\limits_{t_{i}}^{t_{i+1}}
(\lambda(X_{\kappa^a_m(s)} {t_i}+(s- \kappa^a_m(s)))
 - \lambda(X_s))\,ds\right|\,
 |f(X^+_{t_{i}}) - f(X_{t_{i}})|
 \to 0,
\end{eqnarray}
as $\; \Delta\to 0.$

~

Firstly, according to the Lemma \ref{lem3}, there exists a
subsequence $(m'\to\infty)$ such that for almost every $a\in
[0,1]$
\[
\int\limits_0^t |\lambda(X_s) - \lambda(X_{\kappa^a_{m'}(s)})|\,ds \to
0,
\]
and for the whole sequence $(m\to\infty)$,
\begin{equation}\label{la1}
\int\limits_0^1\int\limits_0^t |\lambda(X_s) - \lambda(X_{\kappa^a_m(s)})|\,ds
\,da\to 0.
\end{equation}
Now let us show that
\begin{equation}\label{dif2}
\int\limits_0^1\int\limits_0^t |\lambda(X_{\kappa^a_m(s)} + (s-\kappa^a_m(s))) -
\lambda(X_{\kappa^a_m(s)})| \,ds\,da \to 0, \; m\to\infty.
\end{equation}
The idea is very similar to that of the Lemma \ref{lem3}, however, here a direct reference to this Lemma is questionable and, hence, we need to provide an independent proof. Recall that the trajectory $(X_s)$ may have only finitely many
jumps on $(0,t]$. Also, except for the points of jumps, the
evolution of $(X)$ is deterministic and linear with a constant
positive speed. In other words, integration of the composite
function $\lambda(X_{\cdot})$ turns out to be equivalent to the
integration of $\lambda(\cdot)$ over some interval between two
consequent jumps of the process $(X)$. This allows us to use the
same trick from the Lemma \ref{lem3}, but this time we
approximate the function $\lambda(\cdot)$ rather than the
composite function $\lambda(X_{\cdot})$. Emphasize that this is
possible exactly because of the piecewise linear law of evolution
of $(X_s)$ between the consequent moments of jumps and because
there are only finitely many jumps on each trajectory.

~

\noindent So, let us approximate the function $\lambda(\cdot)$ by
bounded uniformly continuous functions -- say, $\lambda^\epsilon$
-- in the topology of convergence in $L_1[0,T]$ for every $T>0$.
 Note that it is not enough to
approximate $\lambda(\cdot)$ on $(0,t]$. Denote by $(\tau_j, \,
j=0,1,\ldots)$ the sequence of moments of jumps of the trajectory
of $X$; recall that for each $\omega$ there are only finitely many
of them on $(0,t]$.

We have, on each $(\tau_j,\tau_{j+1}]$ and for every $m$,
\begin{eqnarray*}
 \int\limits_0^1\int\limits_{\tau_j}^{\tau_{j+1}} |\lambda^\epsilon
(X_{\kappa^a_m(s)}) - \lambda(X_{\kappa^a_m(s)})|
 \,ds\,da
= \int\limits_{\tau_j}^{\tau_{j+1}} \int\limits_0^1 | \lambda^\epsilon
(X_{\kappa^a_m(s)}) - \lambda(X_{\kappa^a_m(s)})|
 \,da\,ds
 \\\\
\le \int\limits_{\tau_j}^{\tau_{j+1}} 2^m\,\|\lambda^\epsilon - \lambda
\|_{L_1([\tau_j - 2^{-m},\tau_{j+1}])}\,ds
= \|\lambda^\epsilon - \lambda \|_{L_1([\tau_j - 2^{-
m},\tau_{j+1}])} \to 0,\; \epsilon \to 0.
\end{eqnarray*}
Since the number of moments $\tau_j$ on $(0,t]$ is finite for
(almost) all $\omega$, it also follows for each $m$ that
\begin{eqnarray}\label{la2}
 \int\limits_0^1\int\limits_{0}^{t} |\lambda^\epsilon
(X_{\kappa^a_m(s)}) - \lambda(X_{\kappa^a_m(s)})| \,ds\,da \hspace{3.8cm}
 \nonumber \\\\ \nonumber
= \int\limits_0^1\sum_j \int\limits_{\tau_j}^{\tau_{j+1}} | \lambda^\epsilon
(X_{\kappa^a_m(s)}) - \lambda(X_{\kappa^a_m(s)})|\,ds\,da \to 0,\;
 \epsilon \to 0.
\end{eqnarray}

 ~

Now consider the difference
\begin{eqnarray*}
\int\limits_0^1\int\limits_0^t | \lambda^\epsilon(X_{\kappa^a_m(s)} +
(s-\kappa^a_m(s))) - \lambda(X_{\kappa^a_m(s)} +
(s-\kappa^a_m(s)))|\,ds\,da.
\end{eqnarray*}
On each $(\tau_j,\tau_{j+1}]$ and for every $m$, due to
integration over $s$,
\begin{eqnarray*}
\int\limits_0^1\int\limits_{\tau_j}^{\tau_{j+1}} |
\lambda^\epsilon(X_{\kappa^a_m(s)} + (s-\kappa^a_m(s))) -
\lambda(X_{\kappa^a_m(s)} + (s-\kappa^a_m(s)))| \,ds\,da
 \\\\
\le \int\limits_0^1 \|\lambda^\epsilon - \lambda\|_{L_1([\tau_j -
2^{-m},\tau_{j+1}])} \,da = \|\lambda^\epsilon -
\lambda\|_{L_1([\tau_j - 2^{-m}, \tau_{j+1}])} \to 0, \;
\epsilon\to 0.
\end{eqnarray*}
Therefore, also
\begin{eqnarray}\label{la3}
 \left(\int\limits_0^1\int\limits_{0}^{t} | \lambda^\epsilon(X_{\kappa^a_m(s)} +
(s-\kappa^a_m(s))) \right. \hspace{1cm}
\nonumber \\\\ \nonumber
\left.\phantom{\int\limits_0^1}
- \lambda(X_{\kappa^a_m(s)} +
(s-\kappa^a_m(s)))|\,ds\,da \right)\to 0,
\quad  \epsilon\to 0.
\end{eqnarray}
Notice that by virtue of continuity of $\lambda^\epsilon$,
\begin{eqnarray}\label{la4}
\int\limits_0^1\int\limits_0^t | \lambda^\epsilon(X_{\kappa^a_m(s)} +
(s-\kappa^a_m(s))) - \lambda^\epsilon(X_{\kappa^a_m(s)})|\,ds\,da
\to 0, \; m \to \infty,
\end{eqnarray}
for each $\omega$ and $\epsilon$. So, due to (\ref{la2}),
(\ref{la3}) and (\ref{la4}), the left hand side in (\ref{dif2})
goes to zero as $m\to\infty$, as it does not depend on $\epsilon$.
From here it easily follows that there exists a subsequence $(m')$
such that for almost every $a$,
\begin{equation}\label{dif22}
\int\limits_0^t |\lambda(X_{\kappa^a_{m'}(s)} + (s-\kappa^a_{m'}(s))) -
\lambda(X_{\kappa^a_{m'}(s)})| \,ds \to 0, \; m'\to\infty.
\end{equation}
Moreover, this subsequence may be chosen from the earlier fixed
subsequence for which (\ref{la1}) is valid. Hence, we may assume
that for $(m')$ both (\ref{la1}) and (\ref{dif22}) are valid
simultaneously.

~

Overall, by virtue of Lebesgue's dominated convergence, we
conclude that by virtue of (\ref{la1}) and (\ref{dif22}),
convergence (\ref{main_dif1}) holds true with almost any $a\in
[0,1]$, that is, over some subsequence $(m'\to\infty)$,
\begin{eqnarray}\label{main_dif11}
E_X \sum_{i=0}^{2^{m'}} \left|\int\limits_{t_{i}}^{t_{i+1}}
(\lambda(X_{\kappa^a_{m'}(s)})
 - \lambda(X_s))\,ds\right|\,
 |f(X^+_{t_{i}}) - f(X_{t_{i}})|
\to 0. 
\end{eqnarray}
Similarly, without loss of generality, we may  assume that over the same subsequence,
\begin{eqnarray}\label{main_dif111}
E_X \sum_{i=0}^{2^{m'}} \left|\int\limits_{t_{i}}^{t_{i+1}} (h(X_{\kappa^a_{m'}(s)}) - h(X_s))\,ds\right|\,
 |f(X^-_{t_{i}}) - f(X_{t_{i}})|
\to 0.
\end{eqnarray}
From
(\ref{main_dif11})--(\ref{main_dif111}), the statement
(\ref{dynkin1}) (Dynkin's formula) for $f\in C^1_0({\mathcal X})$
follows.

~

\noindent {\bf II.} Now assume $f\in C^1_b({\mathcal X})$. Let us
approximate this function by a uniformly bounded sequence $f^{N}
\in C^1_0({\mathcal X})$, so that
\[
 \|f^N-f\|_{C^1(K)} \to 0, \quad \mbox{for any compact set $K \in {\mathcal X}$}.
\]
Then, by the part I of the proof, we have
\[
E_Xf^N(X_t) - f^N(X) = E_X\int\limits_0^t {\mathcal G}f^N(X_s)\,ds,
\]
and we need to justify the passage to the limit as $N\to\infty$.
This follows from Lebesgue's dominated convergence Theorem.
Indeed, for each $s,t$ and $\omega$,
\[
f^N(X_t) \to f(X_t), \quad {\mathcal G}f^N(X_s) \to {\mathcal
G}f(X_s),
\]
and $f^N(X_t)$ and ${\mathcal G}f^N(X_s) $ are uniformly bounded.
Hence, the statement of the Theorem \ref{thm1} about Dynkin's
formula (\ref{dynkin1}) for any $f\in C^1_b({\mathcal X})$ is
proved.

~

\noindent {\bf III.} Now let us show the formula (\ref{dynkin_t}).
Here we have to consider the difference
\begin{eqnarray*}
 E_X \varphi(t,X_t) - \varphi(0,X) = E_X \sum_{i=0}^{2^m} E_{X_{t_i}}
(\varphi(t_{i+1}, X_{t_{i+1}})-\varphi(t_i,X_{t_{i}})).
\end{eqnarray*}
Consider one term from this sum. This term may be split into five
parts,
\begin{eqnarray*}
E_{X_{t_i}}(\varphi(t_{i+1}, X_{t_{i+1}})-\varphi(t_{i},X_{t_{i}})) \hspace{7cm}
 \\\\
= E_{X_{t_i}}
(\varphi(t_i, X_{t_{i+1}})-\varphi(t_i, X_{t_{i}}))1(\mbox{no
jumps on $(t_{i},t_{i+1}]$})\hspace{2.5cm}
 \\\\
 + E_{X_{t_i}}(\varphi(t_i, X_{t_{i+1}})-\varphi(t_i, X_{t_{i}}))1(\mbox{one
jump up on $(t_{i},t_{i+1}]$})\hspace{1.1cm}
 \\\\
 + E_{X_{t_i}}(\varphi(t_i, X_{t_{i+1}})-\varphi(t_i, X_{t_i}))1(\mbox{one jump
down on $(t_{i},t_{i+1}]$})\hspace{0.6cm}
 \\\\
 + E_{X_{t_i}}(\varphi(t_i, X_{t_{i+1}})-\varphi(t_i, X_{t_{i}}))1(\mbox{at
least, two jumps on $(t_{i},t_{i+1}]$})
 \\\\
+ E_{X_{t_i}}(\varphi(t_{i+1},X_{t_{i+1}})-\varphi(t_i,
X_{t_{i+1}})). \hspace{5.3cm}
\end{eqnarray*}
All terms but the last one are considered similarly to the case
$(f(X))$. The last term, clearly, gives us
\[
E_{X_{t_i}}(\varphi(t_{i+1},X_{t_{i+1}})-\varphi(t_i,
X_{t_{i+1}})) = E_{X_{t_i}}\int\limits_{t_i}^{t_{i+1}}
\frac{\partial}{\partial s} \varphi(s, X_{t_{i+1}})\,ds,
\]
which after summation and taking expectation converges as follows,
\[
E_X\sum_{i=0}^{2^m} E_{X_{t_i}}\int\limits_{t_i}^{t_{i+1}}
\frac{\partial}{\partial s} \varphi(s, X_{t_{i+1}}))\,ds \to E_{X}
\int\limits_{0}^{t} \frac{\partial}{\partial s} \varphi(s, X_{s})\,ds,
\quad m\to\infty,
\]
by Lebesgue's dominated convergence Theorem for any $\varphi\in
C^1_0([0,\infty)\times {\mathcal X})$, as required. Now we may
repeat the arguments about extending the formula from
$C^1_0([0,\infty)\times {\mathcal X})$ to $C^1_b([0,\infty)\times
{\mathcal X})$.

~

\noindent {\bf IV.} It remains to show strong Markov property.
This follows from continuity of the function $E_X f(X_t)$ in $X$
by virtue of the Lemma \ref{Le2} and due to the cadlag property of
the process, see \cite[Theorem 3.4]{Dynkin}. The Theorem
\ref{thm1} is proved.

\section{Proof of Corollary \ref{cor1}}
\noindent Notice that the issue is to extend both
Dynkin's formulae to polynomially growing functions. First of all,
recall that the formulae (\ref{dynkin1}) and (\ref{dynkin_t}) for
$f\in C^1_b( {\mathcal X})$ and $\varphi\in C^1_b([0,\infty)\times
{\mathcal X})$ in martingale language read as follows: the
processes
\begin{equation}\label{mart100}
M_{t}:=f(X_{t}) - f(X) - \int\limits_0^{t} {\mathcal G} f(X_s)\,ds,\; \;
\; t\ge 0,
\end{equation}
and
\begin{equation}\label{mart111}
\tilde M_{t}:=\varphi(t, X_{t}) - \varphi(0,X) - \int\limits_0^{t}
 \left(\frac{\partial}{\partial s}\, \varphi(s,X_s) + {\mathcal
G}\varphi(s,X_s)\right)\,ds,\; \; \; t\ge 0,
\end{equation}
are both martingales (cf. \cite{LSh}). For the reader's
convenience let us show this, say, for $\tilde M$. Indeed,
\[
 E|\tilde M_t| < \infty
\]
since all terms in the right hand side in (\ref{mart111}) are
bounded for $\varphi\in C^1_b$, and
\[
 E_X(\tilde M_{t} - \tilde M_{s}\mid {\mathcal F}^X_s)
 = E_{X_s}(\tilde M_{t-s} - \tilde M_{0}) = E_{X}(\tilde M_{t-s})|_{X=X_s},
\]
due to Markov property. But due to (\ref{dynkin_t}), for {\em
every} initial data $X$,
\[
E_{X}(\tilde M_{t-s})= 0
\]
for every $X$. Hence,
\[
E_{X}(\tilde M_{t-s})|_{X=X_s} = 0.
\]
Therefore, $\tilde M$ is, indeed, a (cadlag) martingale.

~

\noindent Now the task is to extend (\ref{mart111}) to the
function $L_{k,m}$, or, equivalently, to extend (\ref{dynkin_t})
to such function. It is easy to see that $L_{k,m}$ may be
approximated by functions $\varphi^N \in C^1_b([0,\infty)\times
{\mathcal X})$ so that
\[
 \|\varphi^N(\cdot,\cdot) - L_{k,m}(\cdot,\cdot)\|_{C^1([0,N]\times [0,N])} \to 0,
 \quad N\to\infty,
\]
for any $t\ge 0,\, X\in {\mathcal X}$, and
\begin{equation}\label{co1}
 \sup_N \left(|\varphi^N(t,X)| + \left|\frac{\partial}{\partial x} \varphi^N(t,X)\right|
 + \left|\frac{\partial}{\partial t} \varphi^N(t,X)\right| \right) \le C(1+L_{k,m}(t,X)).
\end{equation}
Then it is possible to pass to the limit in all terms of the
equation (\ref{dynkin_t}) written for $\varphi^N$, if the
following a priori bound is established,
\begin{equation}\label{apri}
 \sup_{t\le T} E_{(n_0,x_0)} (n_t+x_t)^m \le C(T,n_0,x_0,m) < \infty,
\end{equation}
for any $m>0$ with some function $C(T,n,x,m)$, for any $T>0$.
Indeed,
\[
\varphi^N(t,X_t) \to \varphi(t,X_t), \quad
\frac{\partial}{\partial s}\varphi^N(s,X_s) \to
\frac{\partial}{\partial s}\varphi(s,X_s), \quad N\to\infty,
\]
and
\[
{\mathcal G}\varphi^N(s,X_s) \to {\mathcal G}\varphi(s,X_s), \quad
N\to\infty.
\]
Since (\ref{apri}) will be established for any $m>0$, we obtain by
Lebesgue's convergence Theorem under the uniform integrability
condition,
\begin{eqnarray*}
\varphi^N(t, X_{t}) - \varphi^N(0,X) - \int\limits_0^{t}
 \left(\frac{\partial}{\partial s}\, \varphi^N(s,X_s) + {\mathcal G}\varphi^N(s,X_s)\right)\,ds
 \\\\
\to \varphi(t, X_{t}) - \varphi(0,X) - \int\limits_0^{t}
 \left(\frac{\partial}{\partial s}\, \varphi(s,X_s) + {\mathcal G}\varphi(s,X_s)\right)\,ds, \;\; N\to\infty.
\end{eqnarray*}
This would also imply that $E_X|\tilde M_t| + E_X|M_t| < \infty$,
for $f=L_{m}$ and $\varphi = L_{k,m}$.

~

\noindent The easiest explanation of (\ref{apri}) is, apparently,
to use the fact that the process $(X_t = (n_t, x_t), \,t\ge 0)$ is
dominated by a similar process, say, $(\bar X_t= (\bar n_t, \bar
x_t), \,t\ge 0)$ {\em without actually serving the customer
at the server} and with a constant arrival
rate $\bar\lambda:= \Lambda \vee \lambda_0$, that is, for each
$\omega$,
\[
\bar x_t\ge x_t, \quad \bar n_t\ge n_t, \quad t\ge 0.
\]
For the process $(\bar X)$ we have,
\[
 \bar x_t = x+t
\]
and
\[
 \bar n_t = n_0 + \xi_t,
\]
where $\xi_t$ has a Poisson distribution with parameter $\bar
\lambda t$. We may imagine a situation as if the current serving has
suddenly pended and the customer remained at the idle server
forever. We have,
\begin{eqnarray*}
 \psi(t,m):=
 E\xi_t^m = \sum_{j=0}^{\infty} j^m\, \frac{(\bar \lambda t)^j}{j!} e^{-\bar \lambda
 t}
 <\infty,
\end{eqnarray*}
so,
\[
 \sup_{t\le T} E_{(n_0,x_0)} (\bar n_t+\bar x_t)^m \le 3^{m-1}(x_0+T)^m +
 3^{m-1} n_0^m + 3^{m-1} \psi(T,m).
\]
Hence, (\ref{apri}) holds true and the statement of the Corollary
\ref{cor1} is proved.

\section{Proof of Theorem \ref{thm2}}
The proof repeats the calculus in \cite{Ve2013_ait} based on
Lyapunov functions
\[
L_{m}(X):= (n+1 + x)^m \quad \mbox{and} \quad L_{m,k}(t,X):=
(1+t)^k(n+1 + x)^m
\]
and on Dynkin's formulae (\ref{dynkin1}) and (\ref{dynkin_t}) due
to the Corollary \ref{cor1}. The news is only a wider class of
intensities, which may be discontinuous, however, this does not
affect the calculus at all once (\ref{dynkin1}) and
(\ref{dynkin_t}) are established. Hence, we drop the details.


\end{document}